\documentclass[11pt]{article}
\usepackage{amssymb,amsmath,amsthm}

\newtheorem{theorem}{Theorem}[section]
\newtheorem{lemma}[theorem]{Lemma}
\newtheorem{proposition}[theorem]{Proposition}
\newtheorem{corollary}[theorem]{Corollary}

\theoremstyle{definition}
\newtheorem{definition}[theorem]{Definition}

\theoremstyle{remark}

\numberwithin{equation}{section}

\newcommand{\e}{\varepsilon}

\newcommand{\ZZ}{{\mathbb Z}}			
\newcommand{\TT}{{\mathbb T}}			

\newcommand{\mcal}{\mathcal}			

\newcommand{\mB}{{\mcal B}}

\newcommand{\mU}{{\mcal U}}

\newcommand{\eset}{\emptyset}			

\newcommand{\Id}{\text{Id}}

\newcommand{\EQ}{\begin{eqnarray*}}
\newcommand{\EQE}{\end{eqnarray*}}

\begin{document}

\title{Sensitive dependence on initial conditions and chaotic group actions}


\author{Fabrizio Polo}


\maketitle

\begin{abstract}
A continuous action of a group $G$ on a compact metric space has
{\it sensitive dependence on initial conditions} if there is a number $\e > 0$ such that for any open set $U$
we can find $g \in G$ such that $g.U$ has diameter greater than $\e.$  We prove that if a $G$ action preserves a probability measure of full support, then the system is either minimal and equicontinuous, or has sensitive dependence on initial conditions.  This generalizes the invertible case of a theorem of Glasner and Weiss.  We prove that when a finitely generated, solvable group, acts and certain cyclic subactions have dense sets of minimal points, then the system has sensitive dependence on initial conditions.  Additionally, we show how to construct examples of non-compact monothetic groups, and transitive, non-minimal, almost-equicontinuous, recurrent, $G$-actions.
\end{abstract}

\section{Introduction}
By a {\it topological dynamical system} we shall mean a continuous action of a topological
group $G$ on a compact metric space $(X,d),$
i.e., a continuous map $\pi: G \times X \to X: (g,x) \mapsto g.x$ with the property that $g.(h.x) = gh.x.$
Our theorems make no reference to the topology on $G.$  So, one may assume $G$ is discrete.
We denote such a system by the pair $(X,G).$  If $T:X \to X$ is continuous, we may also call $(X,T)$ a topological dynamical system.  In this case we are referring to the obvious action $\ZZ.$  Occasionally we refer to theorems in which $T$ is not assumed to be invertible.  It will be made clear when this is the case.

To understand the following property of dynamical systems is the primary goal of this paper:
\begin{definition}
A topological dynamical system $(X,G)$ is said to have {\it sensitive dependence on initial conditions}
(or is {\it sensitive})
if there is some $\epsilon>0$ such that for all $x \in X$ and for every open neighborhood $U$ of $x$ there exists $y \in U$ and $g \in G$ such that $d(g.x, g.y) > \epsilon.$
\end{definition}

It is worth noting that one need not (overtly) mention points in the definition of sensitive dependence on initial
conditions; see the following equivalent definition:
there exists $\e > 0$ such that for all open sets $U, \ \sup_{g \in G} \text{ diam}(g.U) > \e.$

\begin{definition}
A point $x \in X$ is called an {\it equicontinuous point} if for all $\e>0$ there exists $\delta >0$ such that
for any $y \in X, d(x,y) < \delta$ implies $d(g.x,g.y) < \e$ for all $g \in G.$
A point which is not an equicontinuous point will be called a {\it sensitive point.}
\end{definition}

A system $(X,G)$ is said to be {\it equicontinuous} if the set of all maps $\{ x \mapsto g.x : g \in G \}$ is an equicontinuous family. If $(X,G)$ is equicontinuous then it is easily seen that every point is equicontinuous.  Conversely, if every point is equicontinuous then for each $\e>0$ and each point $x$ there exists $\delta_x >0$ such that $g.B_{\delta_x}(x)$ has diameter less than $\e$ for all $g \in G.$  From the balls $B_{\delta_x}(x)$ we can choose a finite subcover $\mU.$ Let $\delta$ be a Lebesgue covering number for $\mU.$  Then, if $d(x,y) < \delta,$ there is some $U \in \mU$ containing $x$ and $y.$  So, $d(g.x,g.y) \leq \text{ diam}( g.U) < \e$ for all $g.$  This proves that every point is equicontinuous if and only if the system is equicontinuous.

Now we define a new metric $d_\infty$ on $X$ by $d_\infty(x,y) = \sup_{g \in g} d(g.x,g.y).$  Consider the identity map $\Id: (X,d) \to (X,d_\infty).$  When $G$ is a monoid, $\Id^{-1}$ is a contraction and hence continuous.  If $(X,G)$ is equicontinuous, it is easy to see that the $\Id$ is a homeomorphism.  Also, if $(X,d_\infty)$ is compact then $\Id$ has compact domain and Hausdorff range and so is a homeomorphism.  As we will see in the proof of Theorem 4, a sensitive point is precisely a discontinuity point for $\Id.$

Logically speaking, the assertion that a system have sensitive dependence on initial conditions is stronger than
the assertion that every point is sensitive (they differ in order of quantifiers.) But, by a category argument (see the excellent article \cite{akinetal}), one can prove
\begin{proposition}
A transitive system is sensitive at each point if and only if it has sensitive dependence on initial conditions.
\end{proposition}

\begin{definition}
A system is called {\it almost equicontinuous} if it has a dense set of equicontinuous points.
\end{definition}

It follows from Corollary 2 in \cite{auslanderyorke} that for a single transitive map, either the system is sensitive, or the set of equicontinuous points is exactly the same as the set of transitive points. This is commonly known as the Auslander-Yorke dichotomy theorem.  It holds in the non-invertible case as well.

The following standard definitions are discussed in much greater detail elsewhere.  For instance, \cite{walters} and \cite{glasner} are good references.  From now on we take $G$ to be a group.
We say $(X,G)$ is {\it topologically transitive} (or {\it transitive} for short)  if for any open $U,V$ there exists $g \in G$ such that $U \cap g.V \neq \eset.$  In our setting, this is equivalent to the existence of a dense set of points $x_0 \in X$ such that $G.x_0$ is dense in $X$ (see Proposition 3.2 from \cite{kontorovich}.)

  A system is called {\it minimal} if every point
has dense orbit.  A point will be called minimal if its orbit closure is a minimal system.  When the acting group is discrete, a point $p \in X$ is called a {\it periodic point} if $G.p$ is finite.  Periodic points are easily seen to be minimal.
Saying that $x$ is a minimal point is equivalent to requiring that for every open neighborhood $U$ of $x,$
$R := \{g \in G : g.x \in U \}$ is {\it left syndetic}. That is, there is a finite subset $F$ of $G$ such that
$FR = G.$

In a group, any set which contains an infinite set of the form $SS^{-1} = \{ st^{-1} : s, t \in S \}$ is called a $\Delta$ set.
A set is called $\Delta^\star$ if it intersects every $\Delta$ set.  A $\Delta^\star$ set which is symmetric (i.e. equal to it's inverse) is always syndetic.
To see this, suppose the set $L \subset G$ is symmetric but not syndetic.  Let $g_0 = 1_G$ and assume we have already chosen distinct elements $g_0, \dots ,g_{n-1}$ with the property that $g_i g_j^{-1} \notin L$ when $i \neq j.$
Since $L$ is not syndetic we can choose $g_n$ such that $\{g_n g_0^{-1}, \dots, g_n g_{n-1}^{-1} \}$ misses $L.$  Since $L$ is symmetric we also know $\{g_0 g_n^{-1}, \dots, g_{n-1}g_n^{-1} \}$ misses $L.$  Inductively, we see $L \cap SS^{-1} = \eset$ where $S = \{g_n \}_n.$  That is, $L$ is not $\Delta^\star.$

A probability measure $\mu$ on the Borel $\sigma$-algebra $\mB = \mB(X)$ is {\it invariant} if $\mu(A) = \mu( g^{-1}.A)$ for all $A \in \mB, g \in G.$  All measures should be assumed to be Borel measures. The measure $\mu$ is said to be {\it ergodic} if any invariant set has measure one or zero. 

In the literature on sensitivity, many results have a common theme: under some additional hypothesis, an almost equicontinuous system must be equicontinuous.  This gives a dichotomy: any system satisfying the additional hypothesis is either equicontinuous or sensitive.  Topological hypotheses are popular (see
\cite{banksetal}, \cite{akinetal}, and \cite{glasnerweiss}.)  Measure theoretic hypotheses also appear (see  \cite{glasnerweiss})
The next two theorems are particularly nice examples.
\begin{theorem} [Glasner, Weiss \cite{glasnerweiss}] \label{glasnerweissThm}
If $(X,T)$ is a transitive topological system equipped with an invariant probability measure $\mu$ of full support
then either $(X,T)$ has sensitive dependence on initial conditions or it is minimal and equicontinuous.
\end{theorem}
\begin{theorem} [Akin, Auslander, Berg \cite{akinetal} ] \label{akinThm}
Let $(X,T)$ be a (possibly non-inervtible) transitive system.  If the set of all minimal points is dense then either the system has sensitive dependence on initial conditions or $X$ is a minimal equicontinuous system.
\end{theorem}
Both of these theorems are generalizations of their predecessor (see \cite{banksetal}):
\begin{theorem} [Banks, Brooks, Cairns, Davis, Stacey \cite{banksetal}] \label{banksThm}
Let $(X,T)$ be a (possibly non-inervtible) transitive system.   If the set of periodic points is dense then either the system has sensitive dependence on initial conditions or $X$ is a finite set.
\end{theorem}
A periodic point is a minimal, so \ref{akinThm} implies \ref{banksThm}. 
Theorem \ref{glasnerweissThm} implies \ref{banksThm} because the existence of a dense set of periodic points allows one to easily construct an invariant measure of full support by adding weighted counting measures on periodic orbits.

In this paper we derive similar results for the actions of more general groups.
Our main results are Theorems \ref{myThm} and \ref{mysolvThm}. Theorem \ref{myThm} generalizes the invertible case of Theorem \ref{glasnerweissThm}.
\begin{theorem} \label{myThm}
Let $(X,G)$ be a transitive system.  If $X$ admits an invariant measure of full support then the system is either minimal and equicontinuous or sensitive.
\end{theorem}

A close examination of the proof in \cite{banksetal} of Theorem \ref{banksThm} reveals that (with obvious modifications) it already works
for any group action.  However, Theorem 2.5 from \cite{akinetal} is more general  so we include it here for completeness.
\begin{theorem} \label{akinThm2}
Let $G$ be a group and $(X,G)$ be transitive system having a dense set of minimal points.
Then, either $(X,G)$ has sensitive dependence on initial conditions or $X$ is minimal, equicontinuous.
\end{theorem}

Heuristically speaking when $G$ is a large group, the requirement that the action have dense periodic or minimal points is quite strong.  Does a result like Theorem \ref{akinThm2} hold if we require only dense minimal points for certain subactions?
Indeed, when $G$ is solvable there is such a theorem.  In order to state it we need another definition.
\begin{definition}
If a finite set $S$ generates a solvable group $G$ we will say $S$ is {\it nice} if
$G^{(n)}$ is generated by $S \cap G^{(n)}$ (where $G = G^{(0)}$ and $G^{(n+1)} = [G^{(n)}, G^{(n)}].)$
\end{definition}

\begin{theorem} \label{mysolvThm}
Let $S$ be a nice generating set for a finitely generated solvable group $G$ which acts transitively on a compact metric space $X.$  Given $s \in S,$ write $\left< s\right>$ for the subgroup generated by $s.$  If the system $(X, \left< s \right>)$ has a dense set of minimal points for each $s,$ then $(X,G)$ is either minimal, equicontinuous or it has sensitive dependence on initial conditions.
\end{theorem}


\section{Almost equicontinuity and transitivity}

In this section we demonstrate the existence of many transitive, almost equicontinuous systems which are not minimal, equicontinuous.
We prove analogues of Theorems 4.2 and 4.6 from \cite{akinglasner}.  The dichotomy type theorems discussed in the introduction use additional hypotheses to eliminate the possibility of almost equicontinuous but not equicontinuous systems.  Without the knowledge that such examples exist, the strength of this kind of theorem is questionable.  

\begin{lemma} \label{inheritanceLemma}
Let $(X,G)$ be a topological dynamical system with an equicontinuous point $x.$ If $y \in X, g_n \in G$ and
$g_n.y \to x$ then $y$ is also equicontinuous and has the same orbit closure as $x.$ 
\end{lemma}

\begin{proof}
Given $\e>0$ there is a neighborhood $U$ of $x$ such that $\text{diam} (g.U) < \e$ for all $g \in G.$  Fix $h$ such that $h.y \in U.$
Then $V:= h^{-1}.U$ is a neighborhood of $y$ with $\text{diam}(g.V) < \e$ for all $g \in G.$  Thus $y$ is an equicontinuity point.
Since $d( g.x, gh^{-1}.y ) < \e$ for all $g \in G,$ we see that the orbit of $x$ is $\e$-dense in the orbit of $y$ and the orbit of $y$
is $\e$-dense in the orbit of $x.$  Thus, taking closures of each orbit yields equal sets.
\end{proof}

\begin{proposition} \label{eqptProp}
Let $(X,G)$ be a system with a transitive point $x.$  The following are equivalent:
\begin{enumerate}
\item $(X,G)$ is almost equicontinuous.
\item $x$ is an equicontinuous point.
\item For all $\e >0$ there exists $\delta>0$ such that for all $g,h \in G,$
$d(h.x,x) < \delta$ implies $d(gh.x,g.x) < \e.$
\item $d$ and $d_\infty$ induce the same topology on the set of transitive points.
\end{enumerate}
\end{proposition}

\begin{proof}
Assume 2.  Then any translate of $x$ is an equicontinuous point.  1 follows.  Now we prove 3 implies 2.  Fix $\e>0$ and choose $\delta>0$ so that $d(h.x,x) < \delta$ implies $d(gh.x, g.x) < \e/2$ for all $g \in G.$  Now fix $g \in G$ and suppose $d(x,y) < \delta/2.$  Choose $h$ to make $h.x$ close enough to $y$ that $d(gh.x, g.y) < \e/2$ and $d(x,h.x) < \delta.$  Then
\[ d(g.x,g.y) \leq d(g.x, gh.x) + d(gh.x, g.y) < \e/2 + \e/2. \]

We will now prove 1 implies 4 and 4 implies 3.  For 1 implies 4, It suffices to show that a sequence $x_n$ of transitive points converges under $d$ to another transitive point $x$, if and only if the same is true under $d_\infty.$  One direction is obvious.  For the other direction, suppose $x_n \to x$ under $d.$  By assumption $x$ is transitive and hence equicontinuous by Lemma \ref{inheritanceLemma}.  For any $\e > 0,$ when $n$ is sufficiently large, we have $d(g.x_n, g.x) < \e$ for all $g \in G.$  That is, $d_\infty(x_n,x) \leq \e.$  Therefore $d_\infty(x_n, x) \to 0,$ as desired.

To see that 4 implies 3, suppose $h_n \in G$ are such that $d(h_n.x, x) \to 0$.  Then $d_\infty(h_n.x,x) \to 0$ as well.  Whence 3.
\end{proof}

\begin{proposition} \label{aeqstructureProp}
Suppose $G$ acts transitively, by isometries on a possibly non-compact metric space $X_0$ and $\iota:X_0 \to X$ is a uniformly continuous metric $G$-compactification.  I.e., $(X,G)$ is a compact metric system and $\iota$ is a uniformly continuous, $G$-equivariant homeomorphic embedding of $X_0,$ onto a dense subset of $X.$  Then $(X,G)$ is an almost equicontinuous, transitive system with $\iota(X_0)$ contained in the transitive points of $X.$  Conversely, every almost equicontinuous, transitive system $(X,G)$ arises in this way: $X_0$ may be taken to be the set of transitive points equipped with the $d_\infty$ metric.
\end{proposition}
\begin{proof}
Since $G$ acts transitively, by isometries on $X_0,$ it acts minimally.  Thus every point $y \in X_0$ has orbit dense in $X_0$ and so $\iota(y) =: x$ has orbit dense in $X.$   We now must show that $x$ is an equicontinuous point.  Fix $\e>0,$ and let $x' := \iota(y').$  Using continuity of $\iota^{-1}$ and uniform continuity of $\iota$ we can choose $\delta > 0$ such that, if $d(x,x') < \delta$ then for all $g \in G,$ $g.y$ and $g.y'$ are sufficiently close that
$d(\iota(g.y), \iota(g.y')) = d(g.x,g,x') < \e.$  In other words, $x$ is an equicontinuous point.

For the converse, let $(X,G)$ be an almost equicontinuous, transitive system and let $X_0$ be the set of all transitive points equipped with the $d_\infty$ metric.  Then the inclusion $\iota:X_0 \to X$ is a contraction and hence uniformly continuous.  By Proposition \ref{eqptProp} part 4, $\iota$ is a homeomorphic embedding.  So, in fact, $\iota$ is a uniformly continuous $G$-compactification, as desired.
\end{proof}

According to Proposition \ref{aeqstructureProp}, constructing almost equicontinuous systems, is equivalent to constructing equicontinuous compactifications of transitive isometric $G$-actions.  One simple way to construct an almost equicontinuous $G$ action is to let $X_0 = G$ with the metric $d(g,h) = 1$ if and only if $g \neq h.$  This is an invariant metric giving the discrete topology.  One could then take the one-point compactification of $G$ and extend the left multiplication action of $G$ by fixing the point at infinity.  The infinite point is sensitive and all other points are equicontinuous and transitive.

From a topological perspective, this example is not very interesting.  Notice that, except for the point at infinity, none of the transitive points are recurrent.  I.e. it is not true that for every transitive point $x,$ every neighborhood $U \owns x,$ and every compact $K \subset G$ we can find $g \notin K$ with $g.x \in U.$

If $G = \ZZ$ we can construct more examples by taking $X = X_0$ to be some compact {\it monothetic} group (a group with a dense cyclic subgroup.)  These are well known and abundant.  They occur precisely as Pontryagin duals of subgroups of the circle equipped with the discrete topology.  Again, these examples are not very dynamically interesting because they are minimal and equicontinuous.

If we take $X_0$ to be a non-compact monothetic group then we can easily pick a metric on it with respect to which $\ZZ$ acts by isometries.  Then any uniform compactification $\iota:X_0 \to X$ gives a transitive almost-equicontinuous system which is not minimal, equicontinuous.  In fact, it is not hard to see that any transitive isometric $\ZZ$-action on a complete metric space $X$ is actually a monothetic group.  This is explored in detail in \cite{akinglasner}.

Now we will show how to construct an example of a $G$-system which is transitive, almost equicontinuous, and recurrent, but not minimal, equicontinuous.  If one analyzes the procedure, we exploit the existence of non-compact monothetic groups.  Some examples of such groups are known (see \cite{rolewicz}) and any of them may be used in our construction.  We will show a different method for constructing such groups.
First we prove some obvious propositions which reduce the problem of defining isometric transitive $G$-actions to defining norms on $G.$

\begin{definition} 
A {\it symmetric norm} on a group $G$ is a function $g \mapsto \| g \| \in [0,\infty)$ such that $\| g \| = 0$ if and only if $g=1,
\| g \| = \| g^{-1} \|,$ and $\|gh \| \leq \|g \| + \| h \|.$  Two norms $\| \cdot \|_i, i=1,2$ are {\it uniformly equivalent} if for all $\e>0$ there is a $\delta>0$ such that for any $i,j$ $\| g \|_i < \delta$ implies $\| g \|_j < \e.$
\end{definition}

\begin{proposition} \label{normcorrespondenceProp}
Transitive, isometric, free $G$-actions on complete, pointed metric spaces are in one-to-one correspondence with symmetric norms on $G.$
\end{proposition}
\begin{proof}
First, suppose $G$ acts transitively, isometrically, and freely on a complete pointed metric space $(X,x).$  Define $\| g \| = d(x, g.x).$  Since the action is free $x \neq g.x$ unless $g = 1.$  Since the action is isometric, $\| g^{-1} \| = d(x, g^{-1}.x) = d(g.x,x) = \| g \|.$  Finally
\[ \| gh \| = d(x, gh.x) \leq d(x, g.x) + d(g.x, gh.x) = d(x,g.x) + d(x,h.x) = \| g \| + \| h \|. \]

Now suppose we have a symmetric norm $\| \cdot \|$ on $G.$  Then we can define a left invariant metric on $G$ by $d(g,h) = \| g^{-1} h \|.$  Let $X$ be the completion of $G$ with respect to this metric and choose $x = 1$ for the base point.  Then $G$ obviously acts on $X$ transitively, isometrically, and freely.
\end{proof}

A point $x$ is said to be {\it recurrent} if for every neighborhood $U \owns x,$ outside every compact subset of $G$ we can find $g$ such that $g.x \in U.$  A $G$-action is said to be {\it recurrent } if for every point is recurrent.

\begin{proposition} \label{normpropertiesProp}
Fix a symmetric norm $\| \cdot \|$ on $G.$
The associated isometric action on the completion $X$ of $G$ is recurrent if and only if there exist $g_n \in G, g_n \to \infty$ such that $\| g_n \| \to 0.$  Furthermore $X$ is compact if and only if for all $\e >0$ the set
$\{ g \in G : \| g \| < \e \}$ is left syndetic.
\end{proposition}
\begin{proof}
Assume $\| g_n \| \to 0.$  Fix any point $x \in X$ and write $x_0 = 1$ for the base point of $X.$  Choose $g \in G$ with $g.x_0 \in B_\e(x).$  Then
\EQ
d(x, g g_n g^{-1}.x) & \leq & d(x, g.x_0) + d(g.x_0, g g_n.x_0) + d(gg_n.x_0, gg_ng^{-1}.x) \\
& < & \e + d(x_0, g_n.x_0) + d( (g g_n g^{-1})g.x_0, (g g_n g^{-1}).x)  \\
& < & \e + d(x_0, g_n.x_0) + \e.
\EQE
But the middle term is equal to $\| g_n \|$ which tends $0.$  So $x$ is recurrent.

For the converse, we use the same argument but reverse the roles of $x$ and $x_0.$  We assume $g_n.x \to x$ and $h_n.x \to x_0.$  We see that we can make $h_n g_n h_n^{-1}$ leave every compact set, while, at the same time
$\| h_n g_n h_n^{-1} \| = d(x_0, h_n g_n h_n^{-1}.x_0) \to 0.$

If $X$ is compact, then $(X,G)$ is a minimal equicontinuous system.  So the set of return times $R$ of $x_0$ to $B_\e(x)$ is left syndetic.  But $R$ is precisely $\{ g \in G : \|g \| < \e \}.$  Conversely, if this set is left syndetic, then we can choose a finite set $F \subset G$ such that $FR = G.$  It follows that for every $\e>0$ and every $x \in X$ we can find $f \in F$ and $g \in R$ such that $d(x,fg.x_0) < \e.$  But
\[ d(x, f.x_0) < d(x, fg.x_0) + d(fg.x_0, f.x_0) < \e + d(g.x_0, x_0) = \e + \| g \| < 2 \e. \]
This proves that the finite set $F.x_0$ is $2\e$-dense.  It follows that $X$ is compact.
\end{proof}

Suppose $\varphi$ is any symmetric nonnegative function on $G$ which takes the value $0$ at $1.$
Define
\[ \| h \| = \inf \{ \varphi(h_1) + \cdots + \varphi(h_n) : h_1 \cdots h_n = h \}. \]
Certainly $\| 1 \| = 0$ and $\| h^{-1} \| = \| h \|.$  Fix $h,h' \in G.$  Given $\e > 0$ we can choose $h_1, \dots, h_{n+k}$ such that $h_1 \cdots h_n = h$ and $h_{n+1} \cdots h_{n+k} = h'$ and
\[ \varphi(h_1) + \cdots + \varphi(h_n) < \| h \| + \e \text{ \ and \ }
\varphi(h_{n+1}) + \cdots + \varphi(h_{n+k}) < \| h' \| + \e. \]
Then $h h' = h_1 \cdots h_{n+k}$ so
\[ \| h h' \| \leq \varphi(h_1) + \cdots \varphi(h_{n+k}) < \| h \| + \| h' \| + 2 \e. \]
It follows that $\| h h' \| \leq \| h \| + \| h' \|.$  Except for the possibility that $\| g \| = 0$ for some $g \neq 1,$ the function $\| \cdot \|$ as defined above is a symmetric norm.

Write $G_\e$ for the subgroup generated by $\{ g \in G : \varphi(g) < \e \}.$
It is easy to see that if, $g \notin G_\e$ then $\| g \| \geq \e.$  So, if for each $g \in G, g \neq 1$ there exists $\e>0$ such that $g \notin G_\e,$ then $\| \cdot \|$ is a symmetric norm.  This condition is far from necessary.  For instance, suppose $G = \ZZ = \left< t \right>$ acts on the circle $X = \TT$ by an irrational rotation.  Define $\| g \|$ to be the distance from $1$ to $g.1.$  Notice that for any $n,$ the subaction generated by $t^n$ is also a minimal.  It follows that for any $\e>0$ we can find $n,k$ such that $t^n$ and $t^{kn+1}$ both have norm less than $\e.$  So $t \in G_\e$ and $G_\e = G.$  Nonetheless, $\| \cdot \|$ is a symmetric norm.

Assume $G$ has an element $t$ of infinite order and define a function $\varphi$ on $G$ as follows.  Let $\varphi(1) = 0.$  For $n \geq 0,$ let $\varphi(t^{n!}) = \varphi(t^{-n!}) = (n + 1)^{-1}.$  Elsewhere, let $\varphi \equiv 1.$  Use $\varphi$ as above to construct $\| \cdot \|.$  Suppose $\| g \| < n^{-1}$ and $g \neq 1.$  Then we can write $g$ as a product of elements of the form $t^{k!}$ where $k \geq n.$  So $g = t^m$ and
$m = a_1 k_1! + a_2 k_2! + \cdots + a_l k_l!$ where the $k_i \geq n$ are distinct, and each $a_i$ is a nonzero integer.  
The assumption on $\| g \|$ tells us that we can do this efficiently so that $\sum_i |a_i| (k_i + 1)^{-1} < n^{-1}.$
Put the $k_i$ in decreasing order and observe that $\sum_{i=2}^l |a_i| < k_1 n^{-1}.$
If $a_1 > k_1+1$ then we could
\begin{enumerate}
\item replace $a_1$ by $a_1-(k_1 + 1)$
\item introduce $k_0 = k_1 + 1$
\item and introduce $a_0 = 1.$
\end{enumerate}
This would give us another way of representing $m$ which reduces the value of $\sum_i |a_i | (k_i + 1)^{-1}.$  Loosely speaking, instead of taking very many large 'steps', we could have taken one even larger step.

A symmetric statement can be made if $a_1 < -k_1-1.$  So, let us assume $|a_1| \leq k_1+1.$  Then
\EQ
|m - a_1 k_1!| &\leq & |a_2| k_2! + \cdots + |a_l| k_l ! \leq (k_1 - 1)! \sum_{i=2}^l |a_i | \\
& \leq &  (k_1 - 1)! \frac{k_1}{n} = \frac{k_1!}{n} .
\EQE
In particular $m$ lies in the interval $[k_1!(a_1-n^{-1}), k_1!(a+n^{-1})].$  Assembling these results over all possible values of $k_1$ and $|a_1| \leq k_1+1$ tells us that $m$ must lie in
\EQ
I := \{ i : \| t^i \| < \frac{1}{n}, i > 0 \}
\subseteq \bigcup_{k \geq n} \bigcup_{a=1}^{k+1} [k!(a-n^{-1}), k!(a+n^{-1})].
\EQE
When $n$ is sufficiently large, $I$ misses any given interval around $0.$  Also $I$ is not syndetic.
Therefore, transferring the statement about exponents to the group itself, we see
\[ \{ g \in G : \| g \| < \frac{1}{n}, g \neq 1 \} = \{ t^i : i \in I \text{ \ or \ } -i \in I \} \]
is not syndetic.  Furthermore, given any element $g \in G, g \neq 1,$ we can choose $n$ large enough that this set does not contain $g.$

Let $X_0$ be the completion of $G$ with respect the metric induced by $\| \cdot \|.$  It is not necessary, but illuminating to observe that this is a disjoint union of non-compact monothetic groups, one for each coset of $\left< t \right>.$  By Propositions \ref{normcorrespondenceProp} and \ref{normpropertiesProp}, $G$ acts transitively, isometrically, freely, and recurrently on the non-compact metric space $X_0.$  Now we may choose any uniformly continuous compactification $\iota: X_0 \to X$ and apply Proposition \ref{aeqstructureProp} to get a transitive, recurrent, non-minimal, almost equicontinuous system $(X,G).$


\section{Invariant measures and ergodicity}

In this section we prove Theorem \ref{myThm}.  First, we need a Lemma (which is interesting in its own right.)

\begin{lemma}
Let $(X,G)$ be a system, and let $\mu$ be an ergodic measure of full support.  Given $A \subseteq X$ of positive measure and $x \in X,$ write $R = R(x,A) = \{ g \in G : g.x \in A \}.$
Then for $\mu$-almost every $x \in X,$  $RR^{-1}$ is $\Delta^\star.$
\end{lemma}
The author would like to thank Vitaly Bergelson for help with this proof.
\begin{proof}
Let $S$ be an infinite subset of $G$ and choose $g=g_S,h=h_S \in S$ such that
$\mu(g^{-1}.A \cap h^{-1}.A) > 0.$  Let $T(g,h)$ be the full measure set
$\bigcup_{k \in G} (kg^{-1}.A \cap kh^{-1}.A).$  If we make such choices for every infinite subset $S,$ countability of $G \times G$ tells us, we have an at most countable collection $T(g_S, h_S)$ of full measure sets.  Let $Y = \bigcap_S T(g_S,h_S).$  Then $\mu(Y) =1.$  Take $x \in Y.$  Then for any infinite $S \subseteq G$ there exist $g=g_S,h=h_S \in S$ and $k \in G$ such that $k.x \in g^{-1}.A \cap h^{-1}.A$.  Therefore $gk.x \in A$ and $hk.x \in A.$  Therefore $RR^{-1} \owns gk(hk)^{-1} = gh^{-1} \in S S^{-1}.$ 
\end{proof}

The proof of Theorem \ref{myThm} was motivated by the techniques in \cite{glasnerweiss}.
\begin{proof}(of Theorem \ref{myThm})
Suppose $(X,G)$ is not sensitive.  That is, it has an equicontinuous point $x.$  Fix $\e>0$ and write $\delta_4 = \e.$  Now choose $\delta_i, i=1,2,3$ such that
\begin{enumerate}
\item If $d(x,y) < 3 \delta_i$ then $d(g.x, g.y) < \delta_{i+1}$ for all $g \in G.$
\item $3 \delta_i < \delta_{i+1}.$
\end{enumerate}

Let $A$ be the $\delta_1$ ball around $x.$ Since $\mu$ has full support we can choose some ergodic component $\nu$ with $\nu(A) > 0.$  Now we can apply Lemma 3 to $\nu$ and $A$ to deduce the existence of a point $y \in A$ with the property that $R := R(y,A)$ satisfies $RR^{-1}$ is $\Delta^\star.$

For $g \in R$ we have
\[ d(g.x,y) \leq d(g.x, g.y) + d(g.y, y) < \delta_2 + \delta_1 < 2 \delta_2.\]
\[ \text{So, \ \ } d(x,g.x) \leq d(x,y) + d(y,g.x) <  \delta_1 + 2 \delta_2 < 3 \delta_2. \]
Taking $h \in R$ we get
\[ d(x,hg^{-1}.x) \leq d(x, h.x) + d(hg^{-1}.(g.x), hg^{-1}.x) <  3 \delta_2 + \delta_3 < \e. \]

We have proven that $RR^{-1} \subseteq R(x,B_\e(x)).$  Since $\e$ was arbitrary, it follows that the set of return times of $x$ to any neighborhood of itself is $\Delta^\star.$  In particular, $R(x,B_{\delta_3}(x))$ is $\Delta^\star.$
Given $g \in R(x,B_{\delta_3}(x)),$ we know $d(x, g^{-1}.x) = d(g^{-1}.(g.x), g^{-1}.x) < \e.$  It follows that
$R(x,B_{\delta_3}(x)) \cup R(x,B_{\delta_3}(x))^{-1} \subseteq R(x,B_\e(x)).$  In other words, $R(x,B_\e(x))$ contains a symmetric $\Delta^\star$ and so must be (left) syndetic (as explained in the introduction.)
But, this is equivalent to the assertion that $x$ have minimal orbit closure (see for instance \cite{glasner}.).  Since $x$ is transitive, we conclude that $(X,G)$ is minimal.  By Lemma \ref{inheritanceLemma}), $(X,G$) is equicontinuous.
\end{proof}

In the case that $\mu$ is an ergodic measure, there is an alternate proof of Theorem \ref{myThm} that relies on Proposition \ref{aeqstructureProp}.  The method is completely different and requires a simple lemma.  First define the function
\[ s(x) = \inf_{\e > 0} \sup \{ \text{ diam}(g.B_\e(x)) : g \in G\} = \inf \{ d_\infty \text{-diam}( U ) : x \in U \text{ open} \} \]

We call $s(x)$ the {\it sensitivity constant} at $x.$  Notice that $x$ is a sensitive point if and only if $s(x) > 0.$

\begin{lemma} \label{sLemma}
The function $s$ defined above is upper semi-continuous and hence measurable.  If $(X,G)$ admits an ergodic measure $\mu$ of full support for which the set of all sensitive points has positive measure then the system has sensitive dependence on initial conditions.
\end{lemma}
\begin{proof}
Suppose $U$ is a neighborhood of $x$ such that $d_\infty \text{-diam}(U) < s(x) + \e.$  Then for any $y \in U,$
$s(y) \leq d_\infty \text{-diam}(U) < s(x) + \e.$  Therefore $s$ is upper semi-continuous.  Measurability follows.

Since $s$ is an invariant function it is constant $\mu$-almost everywhere.  Since it was assumed to take positive values on a set of positive measure, it must be equal to some $c>0,$ $\mu$-almost everywhere.  Thus $s^{-1}(c)$ is a dense set.  Let $U$ be any open subset of $X.$  Then $U$ contains an element of $s^{-1}(c).$  By definition of $s$ we can find $g \in G$ such that $\text{diam} ( g.U ) > c/2.$  So, $(X,G)$ is sensitive.
\end{proof}

Let $(X,G)$ be an almost-equicontinuous system and let $\mu$ be an ergodic probability measure of full support.  Then the set $X_0$ of equicontinuity points is an invariant set and so must have measure $1$ or $0.$  If it has measure $0$ then almost every point is sensitive.  By Lemma \ref{sLemma}, $(X,G)$ is sensitive.  If $X_0$ has measure $1$ then by Proposition \ref{aeqstructureProp} we can think of $\mu$ as a Borel measure on $(X_0, d_\infty).$

If $(X_0, d_\infty)$ is not compact then for some $\e>0$ we can choose a sequence $x_1, x_2, \dots \in X_0$ with $d_\infty(x_i, x_j) > \e$ when $i \neq j.$  Cover $X_0$ by countably many balls of $d_\infty$-radius $\e/4.$  One of them must have positive measure.  Call it $B$ and choose $g_n \in G$ such that $x_n \in g_n.B.$  Since $G$ acts on $(X_0, d_\infty)$ by isometries, the balls $g_n.B$ are disjoint.  This is ludicrous, since they must each have positive measure.

So, $(X_0, d_\infty)$ must be compact.  Continuity of $(X_0,d_\infty) \to (X,d)$ tells us $X_0$ is compact as a subspace of $X.$  Density tells us $X_0 = X.$  The identity is then a homeomorphism and $(X,G)$ is isomorphic to an isometric system.  I.e. $(X,G)$ is minimal and equicontinuous.


\section{Periodic points and minimal subsystems}

In \cite{devaney},  Devaney suggests three properties which define the essence of chaos.  According to him,
a dynamical system (i.e. a continuous map $T:X \to X$) should be called chaotic if it is
\begin{enumerate}
\item topologically transitive,
\item has a dense set of periodic points,
\item has sensitive dependence on initial conditions.
\end{enumerate}

It was first observed in \cite{banksetal} that these requirements are not independent.  In fact, they prove that the first
two conditions imply the third (this is the content of Theorem \ref{banksThm}.)  In \cite{glasnerweiss}, Glasner and Weiss derive this as a corollary of Theorem \ref{myThm}.  They also produce a remarkably simple direct proof (their Corollary 1.4.)  Unfortunately, it is unclear to this author how to adapt the second argument to the case of a non-abelian acting group.

Now we set out to prove Theorem \ref{mysolvThm}.

\begin{lemma} \label{niceLemma}
Assume $S$ is a nice generating set for the solvable group $G.$  Then there is some bound $N$
so that any $g \in G$ can be written in the form
$g = s_{i_1}^{e_1} s_{i_2}^{e_2} \cdots s_{i_N}^{e_N} $ where each $s_i \in S, e_i \in \ZZ.$
\end{lemma}

\begin{proof}
If $G$ is abelian this is obvious, so assume $G$ has higher solvability degree.
Enumerate $S = \{s_1, s_2, \dots s_n \}.$  We know we can
write $g = s_{i_1}^{e_1} s_{i_2}^{e_2} \cdots s_{i_k}^{e_k} $ for some $k.$  Rearrange this word
by collecting like terms and introducing commutators where necessary.
That is, rearrange the word into the form:
\[ g = s_1^{f_1} c_1 s_2^{f_2} c_2 s_3^{f_3} c_3 \cdots c_{n-1} s_n^{f_n} c_n \]
where $c_j \in [G,G].$  By induction on degree of solvability, there is some constant $M$
such that each $c_i$ contributes at most $M$ terms of the form $s_a^b.$  So, $g$ can be written
with no more than $nM + n = |S|(M +1)$ terms.
\end{proof}

\begin{lemma} \label{inheritance2Lemma}
Suppose $x$ is an equicontinuous point in a topological dynamical system $(X,G).$  If $x$ is a limit
of minimal points $x_n$ then $x$ is minimal.
\end{lemma}
This lemma appears in \cite{akinetal}.
\begin{proof}
Fix $\e>0$ and choose $\delta>0$ smaller than $\e/2$
such that $g.B_{\delta}(x) \subseteq B_{\e/2}(g.x).$  Choose $y := y_n \in B_{\delta}(x)$ and
let $R = \{ g \in G : g.y \in B_{\e/2}(x).$  Since $y_n$ is minimal, $R$ is syndetic.  If $g \in R$ then
$d(g.x, x) \leq d(g.x, g.y) + d(g.y, x) \leq \e/2 + \e/2.$  So
$R \subseteq \{ g \in G : g.x \in B_{\e}(x) \},$ which proves the latter set is left syndetic.
Thus, $x$ is minimal.
\end{proof}

\begin{proof}(of Theorem \ref{mysolvThm})
Suppose the system is not sensitive.  Then it has an equicontinuous point.  By Theorem \ref{akinThm}, any transitive point must be equicontinuous.  Let $x$ be such a  point.    Fix $y \in X = \overline{ G.x }.$
Choose a sequence $g_n \in G$ such that $g_n.x \to y.$  Lemma \ref{niceLemma} tells us that each $g_n$ can be written with at most $N$ terms of the form $s^e.$  By passing to a subsequence we can find a sequence $s_1, s_2, \dots, s_k \in S$ such that each $g_n$ is of the form $g_n = s_k^{e_{n,k}} s_{k-1}^{e_{n,k-1}} \cdots s_1^{e_{n,1}}$
(we allow the possibility that $s_i = s_j$ for unequal $i,j.$)

Let $X_1$ be the orbit closure of $x$ under $\left< s_1 \right>.$
Pass to a subsequence along which $s_1^{e_{n,1}}.x$ converges to some point $x_1 \in X_1.$
The points which are minimal under the action of $\left< s_1 \right>$ are dense by assumption.  So, by Lemma \ref{inheritance2Lemma}, $x$ is
also minimal under this action.  Thus we can find powers of $s_1$ which move $x_1$ arbitrarily close to $x.$  Applying Lemma \ref{inheritanceLemma} we see that $x_1$ is another transitive equicontinuous point.
Since $x$ is an equicontinuous point, as $d(s_1^{e_{n,1}}.x, x_1) \to 0$ we also have
\[ d(  s_k^{e_{n,k}} \cdots s_{2}^{e_{n,2}}.( s_1^{e_{n,1}}.x ),     s_k^{e_{n,k}} \cdots s_{2}^{e_{n,2}}.x_1) \to 0. \]
This proves that $s_k^{e_{n,k}} \cdots s_{2}^{e_{n,2}}.x_1 \to y.$

Repeat this process: let $X_2$ be the orbit closure of $x_1$ under $\left< s_2 \right>.$
Pass to a subsequence along which $s_2^{e_{n,2}}.x_1$ converges to some point $x_2 \in X_2.$
Argue as above to conclude that $X_2$ is minimal under the action of $\left< s_2 \right>,$
$x_2$ is transitive equicontinuous, and satisfies
$s_k^{e_{n,k}} \cdots s_{3}^{e_{n,3}}.x_2 \to y.$

Continuing in this way we get a transitive equicontinuous point $x_{k-1},$ whose orbit closure $X_k$ under
the action of $s_k$ is minimal and contains $y.$  Then we can find some powers of $s_k$ which move $y$
as close as we like to $x_{k-1}.$  By Lemma \ref{inheritanceLemma}, $y$ is a transitive, equicontinuous point.

We have shown that every point is transitive and equicontinuous.   Equivalently, $(X,G)$ is a minimal equicontinuous system.
\end{proof}

\begin{corollary}
Let $G,S$ be as in Theorem \ref{mysolvThm}  and assume the set of periodic points for $(X, \left< s \right>)$ is dense for each $s \in S$ (where $\left< s \right>$ is the group generated by $s.$)  If $(X,G)$ is not sensitive then $X$ must be a finite set on which $G$ acts transitively.
\end{corollary}

It would be advantageous to drop the condition that $S$ be nice in Theorem \ref{mysolvThm}.  Unfortunately, the following example
demonstrates that this condition (or something like it) is unavoidable.
Let $G$ be the solvable group $\ZZ \rtimes \ZZ/2\ZZ$ and let $a$ and $b$ be the generators of the factors $\ZZ$ and $\ZZ/2\ZZ$ respectively.  Then $S := \{a, ab \}$ is another generating set for $G.$
Let $X$ be the one point compactification of $G$ and extend the left multiplication action of $G$ on itself to an action on $X$ by fixing the point at infinity.  The system $(X,G)$ is transitive but not minimal.

Each $s \in S$ has order two. So, every point is minimal for the system $(X,\left< s\right>).$  However,
$(X,G)$ does not have sensitive dependence on initial conditions.  In fact, the only sensitive point is the point at infinity.  All other points are isolated and therefore equicontinuous.  This is not a counterexample to Theorem \ref{mysolvThm} because $S$ is not nice.  Suppose we add more generators to $S$ to make it nice.  For instance we could take $S = \{a,ab, [a,ab] \}.$  Then the hypotheses of the theorem are not met: for any $x \in X, \lim_{|n| \to \infty} [a,ab]^n.x = \infty,$ which proves the only minimal point of $(X,\left<[a,ab]\right> )$ is $\infty$ (clearly not dense.)
\newline \newline
Acknowledgements:
The author would like to thank Vitaly Bergelson for proposing this avenue of research, Cory Christopherson for proof-reading an early version, and Michael Hochman, who made me aware of \cite{akinetal}.  The author would especially like to thank Rafal Pikula for thoroughly checking the final draft.

\
\bibliographystyle{amsplain}

\end{document}